\documentclass[11pt,psamsfonts]{amsart}

\usepackage{amssymb}

\numberwithin{equation}{section}

\begin{document}

\title{$L_1$ Compactness of bounded $BV$ Sets}

\author{Isidore Fleischer}
\address{Centre de recherches math\'ematiques, C.P. 6128,
succursale Centre-ville, Montr\'eal, Qc, H3C 3J7, Canada}


\begin{abstract}
Functions, uniformly bounded in $BV$ norm in some bounded
open set $U$ in $R^n$, are compact in $L_1(U)$. This result
is known when $U$ has Lipschitz boundary [EG Th. 4 p.
176], [G 1.19 Th. p. 17], [Z 5.34 Cor. p. 227]; the proof
for general $U$ here, after identifying the operator
theoretic definition of bounded $BV$ norm with that
of the Tonelli variation, appeals to the standard
compactness criterion in $L_1$ [DS 21 TH. p. 301]
[Y, p. 275] (For completeness, these two auxiliary
results are also presented).
\end{abstract}
\maketitle

\section*{Proof}

For the purpose of establishing the compactness criterion,
it will be necessary to bound integrals of the form
$\int^1_0V(x+h) - V(x)\, dx$ for monotone non-decreasing
$V$ on $[0, 1]$, constant on $[1, \infty)$, for $h > 0$.
Since
$\int^1_0V(x+h)\, dx = \int^{1+h}_hV(x)\, dx$, the
difference integrates to $\int^{1+h}_1V(x)\, dx -
\int^h_0V(x)\, dx$, bounded by $2h[V(1)-V(0)]$. It thus
appears that the integral of the differences goes to 0 with
$h$, uniformly in the rise $V(1) - V(0)$ of $V$. If $V$ is
the variation of a $BV$ function $F$, then the same bound
holds for $\int|F(x+h)-F(x)|\, dx$. If the $V$ are functions
of other variables, one can also obtain this uniformity for
the integrals extended over bounded (in measure) subsets of
these.

\section*{The Compactness Criterion}

A bounded subset $\mathcal K$ of $L_1(U)$ ($U$ open
bounded), satisfying\linebreak $\lim_{t\to
0}\int|f({s+t})-f(s)|ds = 0$ uniformly for $f \in \mathcal
K$ is precompact.

The hypothesis entails (by $s \to s +t_0$) the formally
stronger $\int|f(s+t) - f(s+t_0)|\, ds \to 0$ for every
$t_0$ and uniformly in $t_0$ as well as in $\mathcal K$.

$(M_rf)s := \text{-}\mskip-12mu\int_{B_r} f(s+t)\, dt$
($\text{-}\mskip-12mu\int_E := (1/|E|)\int$, average of
integrand over measurable
$E$; $B_r :=$ centered ball of radius $r$) is uniformly
equicontinuous in $\mathcal K$ and $r \geq r_0 > 0$ and has
the same $L_1$ norm as $f$ -- hence, this set of $M_rf$ is
precompact (in $C$, by Ascoli, hence) in $L_1$.

\noindent $\|M_rf-f\| =
\int|\,\text{-}\mskip-12mu\int_{B_r}f(s+t)-f(s)dt|ds
\leq f_{B_r}\int|f(s+t)-f(s)|\, ds\, dt\leq \sup_{t\in B_r}
\int|f(s+t)-f(s)|\, ds \to 0$
as $r \downarrow 0$, uniformly in $\mathcal K$. It thus
suffices to verify that $\{M_rf: r > 0, f \in \mathcal
K\}$ is precompact and for this, that every
sequence with $r \to 0$ has a subsequence every two terms of
which are $< 1/n$ apart. For any such sequence, find an
$r_0$ below which $M_rf$ is  within $1/2n$ of $f$ uniformly
in $\mathcal K$ and use the precompactness of the $M_{r_0}f$
to extract such a subsequence. By diagonalizing, this will
produce a Cauchy subsequence.

\section*{$BV$ Varieties}

The (more precisely, ``A'') \emph{variation} of an $n$-
argument real-valued locally integrable function $F$ in an
open bounded set $U$ is defined to be the sup of
$\int_U F\bigtriangledown\cdot\varphi\, dx$ over the $C'\,
\varphi$ to $[-1+1]$ with compact support in $U$.

The integral is a sum of $n$ terms $\int_U
F\partial\varphi_{\imath}/\partial x_{\imath}\,
dx_{\imath}\, d\tilde x_{\imath}$ where $d\tilde
x_{\imath}$ is the element of $(n-1)$-volume in the
co-ordinate hyperplane orthogonal to the $x_{\imath}$- axis
and $\varphi_{\imath}$ is the $\imath^{th}$ component of
$\varphi$. Since $\varphi$ could have just this component
different from 0, the sup of the sum dominates the sup of
each term; and the sum of these non-negative sups dominates
the sup of their sum. Hence bounded variation comes to
finiteness of the $n$ sups of $\int_U
F\partial\varphi_{\imath}/\partial x_{\imath}\,
dx_{\imath}\, d\tilde x_{\imath} =
\int_U-\varphi_{\imath}\, d_{\imath}F\, d\tilde x_{\imath}$
where $d_{\imath}F$ is the distribution in $x_{\imath}$
obtained from $F$ by fixing the other co-ordinates; the
minus sign can be dropped since the $\varphi$'s are closed
for negation. The value of the integral is unchanged by
making $F$ right continuous in $x_{\imath}$, so assume it
so. The sup is dominated by the ``Tonelli Variation'' $\int
V_x(F)\, d\tilde x$, the integrand being the classical
variation of right continuous $F$ with all the variables
other than $x$ fixed.\footnote{The integrand is positive
measurable as limit of a sequence of such; a priori the
integral could be infinite.} This quantity is actually
attained by the sup as we now show.

It would suffice to attain an arbitrary approximating
sum whose terms are measures of finitely many disjoint
$\tilde x$-measurable subsets with values of the integrand
at points in the sets as coefficients: this integrand is
itself a sup of approximating sums: $\Sigma|F(t_{1+1}) -
F(t_{\imath})|$ with $t_{\imath}$ points of continuity of
$F$.

It is easy to find a sequence of $[-1, +1]$-valued
$C'$ functions on a real interval $[s_{\imath} t]$
vanishing, along with their derivatives, at $s$ and $t$
which increase pointwise to $1_{(s_{\imath} t)}$; the
integral of this sequence against $BV\, dF$ converges to
$F(t^-) - F(s^+)$. Glueing together translates of this
sequence or its negative yields a $C'$ sequence to $[-1,
+1]$ whose integral against $dF$ converges to
$\Sigma|F(t_{1+1}) - F(t_{\imath})|$ provided the partition
points $t_{\imath}$ are restricted to points of continuity
of $F$.

One can also find a sequence of $[-1, +1]$- valued $C'$
functions on $\widetilde X$ which converges a.e. to
the characteristic function of a given measurable set:
represent the set and its complement a.e. as a union of a
sequence of closed sets; for like-indexed pairs find (by
Urysohn) a continuous function 1 on the contained closed
set and 0 on the one disjoint; and approximate the
resulting sequence (by Weierstrass-Stone) by a $C'$
sequence with the same pointwise limit.
The $\tilde x$-integral of these functions converges to
the measure of the set (which the function characterizes)
and so their pairwise product with the real ($x$-)
variable sequence previously constructed yields a $C'$
sequence on $\widetilde X \times X$ to
$[-1, +1]$ whose $\int\varphi\, dF\, dx$ converges to
$\Sigma|F(t_{i+\imath}) - F(t_i)|$ multiplied by the
measure of the $\tilde x$-subset, which is one term of
the sought for approximating sum; adding these finitely
many sequences pointwise yields a sequence converging to
the full sum.

\vspace{0.5cm}
\noindent P.S. The same result is valid  with the usual
(Vitali) notion of $BV$: Over a bounded
domain, a uniformly $BV$, $L_1$-bounded set of
functions is uniformly bounded (for if the values at some
point were unbounded then by $BV$ uniformity so would be
their lower bounds uniformly over the domain; and then
again by uniform
$BV$, the pointwise bounded set would be uniformly bounded
over the domain). By Helly's Selection Theorem, some
subsequence converges pointwise, whence by the uniform
boundedness, in $L_1$.  

\bibliographystyle{amsalpha}

\begin{thebibliography}{WZ}

\bibitem[DS]{DS} N.~Dunford and J. Schwartz, \emph{Linear
Operators}, Pt 1. Interscience.

\bibitem [EG] {EG} L. C. Evans and R. F., Gariepy,
\emph{Measure theory and fine properties of functions}, CRC
Press, 1992.

\bibitem [G] {G} E. Giusti, \emph{Minimal Surfaces and
Functions of Bounded Variation}, Birkh\"{a}user, 1984.

\bibitem[Y] {Y} K.~Yosida, \emph{Functional Analysis},
Springer Verlag, 1968.

\bibitem [Z] {Z} W. P. Ziemer, \emph{Weakly differentiable
functions}, Springer-Verlag 1989.
\end{thebibliography}

\end{document}